\documentclass[10pt,english]{article}
\usepackage{subfigure}
\usepackage{graphicx}
\usepackage{float}
\usepackage[T1]{fontenc}
\usepackage[latin9]{inputenc}
\usepackage[margin=1.25in]{geometry}
\usepackage{float}
\usepackage{amsmath}
\usepackage{amsbsy}
\usepackage{setspace}
\usepackage{amssymb}
\usepackage{esint}
\usepackage{cite}
\newfloat{algorithm}{tbp}{loa}
\floatname{algorithm}{Algorithm}
\usepackage{hyperref}
\usepackage{listings}

\newlength\myindent
\makeatletter

\floatstyle{ruled}
\newfloat{algorithm}{tbp}{loa}
\floatname{algorithm}{Algorithm}

\usepackage{babel}

\begin{document}

\title{Heuristic Optimal Transport in Branching Networks}

\author{M. Andrecut}

\date{November 11, 2023}

\maketitle
{

\centering Unlimited Analytics Inc.

\centering Calgary, Alberta, Canada

\centering mircea.andrecut@gmail.com

} 
\begin{abstract}
Optimal transport aims to learn a mapping of sources to targets by minimizing the cost, which is typically defined as a function of distance. 
The solution to this problem consists of straight line segments optimally connecting sources to targets, and it does not exhibit branching. 
These optimal solutions are in stark contrast with both natural, and man-made transportation networks, where branching structures are prevalent. 
Here we discuss a fast heuristic branching method for optimal transport in networks. We also provide several numerical applications to 
synthetic examples, a simplified cardiovascular network, and the "Santa Claus" distribution network which includes 141,182 cities around the world, 
with known location and population. 
\end{abstract}

\section{Introduction}

The Optimal Transport (OT) problem was first introduced by Monge \cite{key-1}, and recently it has gained an extensive attention due to its applications in machine learning, biology and economics \cite{key-2}, \cite{key-3}, \cite{key-4}. 
Given two probability measures $\mu$ and $\nu$ over the metric spaces $\mathcal{X}$ and $\mathcal{Y}$, and a cost function $c:\mathcal{X}\times \mathcal{Y} \rightarrow \mathbb{R_{+}}$ 
of transporting a unit mass from $x \in \mathcal{X}$ to $y \in \mathcal{Y}$, the modern formulation of OT seeks to find the map $\phi:\mathcal{X} \rightarrow \mathcal{Y}$ satisfying:
\begin{equation}
\text{arg}\inf_\phi \left\lbrace \int_\mathcal{X} c(x,\phi(x))d\mu(x) \mid \phi_\#\mu = \nu \right\rbrace,
\end{equation}
where $\phi_\#$ denotes the push-forward map of $\mu$ by $\phi$. An important result by Brenier \cite{key-5} shows that if the cost function is the quadratic distance $c(x,y)=\Vert x -y \Vert^2$, 
and at least one of $\mu$ and $\nu$ has a density with respect to the Lebesgue measure, then the solution of the OT problem exists and is unique. Therefore, if 
these conditions are not met, then the problem may not always have a solution. However, the convex relaxation introduced by Kantorovich \cite{key-6}:
\begin{equation}
\text{arg}\inf_\gamma \left\lbrace \int_{\mathcal{X}\times \mathcal{Y}} c(x,y)d\gamma(x,y) \mid \gamma  \in \Gamma(\mu, \nu) \right\rbrace,
\end{equation}
is guaranteed to have a solution. Here, $\Gamma(\mu, \nu)$ is the set of "transportation plans" defined by the joint distributions with marginals $\mu$ and $\nu$. 

The Kantorovich formulation is more realistic, and it provides the support for modern transport and resource allocation problems. However, the solution to this problem 
consists of straight line segments (or geodesics) optimally connecting sources to targets, which is in stark contrast with both natural, and man made transportation networks, 
where branching structures are prevalent. The reason behind this is that in the Monge-Kantorovich formulation the problem only considers sources and targets, 
but in reality it may be more optimal to also use intermediary branching points, such that the items are first transported together to a branching point and then they are 
distributed to targets. For example, in the case of one source and two targets the OT solution is a V shaped path, however a Y shaped branching solution may be better. 
Many natural and man-made transportation networks, such as veins on a leaf, blood circulation, water and gas supply, or the post office mailing system, actually require branching 
in order to perform optimally. Thus, an important objective is also to find new sites (branching points) relative to the set of already 
existing targets in order to further improve resource allocation in complex systems.

Here we provide a fast heuristic method for Branching Optimal Transport (BOT). The method is inspired by fluid networks, and it is 
based on a local optimization algorithm, defined as a combination of greedy and tabu search, which can provide fast practical solutions to larger scale transport problems. 
We assume that the sources and the targets are discrete, but the fluid to be transported is continuously distributed and flowing through the network, 
and therefore it is not localized in point like objects. 
This means that given the required flow distribution we are interested in finding the optimal branching structure of the network connecting the sources and the targets. 
We show that with these assumptions, a Y branching solution is interpolating between a T like branching solution, 
and the V like no branching corresponding to the OT solution. Once the branching network is constructed by the algorithm the "fluid" constraint can be 
removed, and the transportation network can be used for both point like objects and fluids. 
We also discuss several numerical applications, including synthetic examples, a simplified cardiovascular network, and the "Santa Claus" distribution network which includes 141,182 cities around the world, 
with known location and population.

\section{Discrete optimal transport}
In many practical applications the source and target distributions $\mu$ and $\nu$ are only given by discrete samples $X=\{x_i\}_{i=1}^{N_x}$, $Y=\{y_j\}_{j=1}^{N_y}$:
\begin{equation}
\mu = \sum_{i=1}^{N_x} p_i \delta(x_i), \quad \nu = \sum_{j=1}^{N_y} q_j \delta(y_j),
\end{equation}
where $\delta(x)$ is the Dirac function, and $p_i$, $q_j$ are the probabilities associated with the samples:
\begin{equation}
\sum_{i=1}^{N_x} p_i = 1, \quad \sum_{j=1}^{N_y} q_j = 1, \quad p_i,q_j > 0.
\end{equation}
In this case the transportation cost and the plan are the matrices $\gamma = [\gamma_{ij}]_{N_x \times Ny}$, $c = [c_{ij}]_{N_x \times Ny}$, and the total cost is given by their Frobenius product:
\begin{equation}
C = \langle \gamma,c \rangle= \sum_{i=1}^{N_x} \sum_{i=j}^{N_y} \gamma_{ij} c_{ij}.
\end{equation}
Therefore, the discrete version of the OT problem aims to find the optimal transportation plan by solving:
\begin{equation}
\gamma^* = \min_{\gamma} \langle \gamma,c \rangle,
\end{equation}
with the constraints defined above \cite{key-2}, \cite{key-3}.
This is a linear programming problem, however in many practical cases it could be very large. 
In such cases, it has been shown that it is more efficient to solve the entropic regularized problem:
\begin{equation}
\gamma^* = \min_{\gamma} \langle \gamma,c \rangle - \lambda H(\gamma), \quad \lambda > 0.
\end{equation}
This regularized problem can be solved much faster using the Sinkhorn method, which is an iterative matrix-scaling algorithm \cite{key-7}. 

While this approach does find the optimal transport plan, it does not solve, or explain, the branching problem. The branched transportation network 
can be defined as a directed, edge-weighted acyclic graph $G(V,E)$, with the edges $E\subset V\times V$ connecting the nodes $V=X\cup Y \cup Z=\{v_n\}$, where $X=\{x_i\}$ are the sources, 
$Y=\{y_j\}$ the targets, and $Z=\{z_k\}$ the intermediary branching points \cite{key-8}. 
Assuming that $m_{ij}$ is the mass transported on the edge $(i,j)\in E$, the following balancing equations can be written at each node:
\begin{itemize}
\item Supply for sources, $i=1,...,N_x$:
\begin{equation}
p_i = \sum_k m_{ik} - \sum_k m_{ki}
\end{equation}
\item Demand for targets, $j=1,...,N_x$:
\begin{equation}
q_j = \sum_k m_{kj} - \sum_k m_{jk}
\end{equation}
\item Conservation for branching points, $\ell=1,...,N_z$:
\begin{equation}
\sum_k m_{k\ell} = \sum_k m_{\ell k}
\end{equation}
\end{itemize}
The optimal branching points $Z$ and the corresponding mass distribution $m=\{m_{ij}\}$ must be the solution of the following BOT optimization problem:
\begin{equation}
\text{arg}\min_{Z,m} \sum_{(i,j)\in E}  m_{ij}^\alpha \Vert v_i - v_j \Vert, \quad \alpha \in [0,1].
\end{equation}
For $\alpha=1$ the BOT problem is equivalent to the discrete OT problem, while for $\alpha = 0$ it is equivalent to the Euclidean Steiner problem \cite{key-9}, which requires the 
shortest network interconnecting all the nodes independently of the mass distributions. For intermediary values of $\alpha \in (0,1)$ the BOT problem interpolates 
between these two optimization problems, and it has been shown to be NP-complete \cite{key-10}.

\section{Local optimization}

Most of the current research on the BOT problem is based on the sub-additivity rule $m\rightarrow m^\alpha$ described in the previous section, since it captures 
the increased efficiency of aggregate transportation:
\begin{equation}
(m_{ij} + m_{ik})^\alpha < m_{ij}^\alpha + m_{ik}^\alpha, \quad \alpha \in [0,1),
\end{equation}
and it has been shown that the solutions indeed exhibit a branching structure \cite{key-11}, \cite{key-12}.  
This approach has generated several more or less complicated methods for locally computing the branching angles, having different requirements derived from various optimal principles (see \cite{key-11}, \cite{key-12} 
and the references within). 
Here we consider a more simplified approach.

Given the flow from a node at location $v_k$ to the nodes $v_i$ and $v_j$, the local optimization problem we would like to solve is to find the 
location $z_b$ of a branching point $b$ such that the following local optimal condition is satisfied:
\begin{equation}
z_b = \min_{z} \left\lbrace  F(z) = m_{kb}\Vert v_k - z \Vert + m_{bi}\Vert v_j - z \Vert + m_{bj}\Vert v_j - z \Vert \right\rbrace. 
\end{equation}
We approximate the mass of each branch by the mass of a cylinder with density $\rho$, sectional area $s$ and length $\ell$: $m(\rho,s,\ell)=\rho \ell s$. 
Also we assume that the sub-additivity rule applies only to the sectional area: $f(s) \rightarrow s^\alpha$.
After several simplification steps we obtain the following quadratic functional: 
\begin{equation}
F(z) = s_{kb}^\alpha\Vert v_k - z \Vert^2 + s_{bi}^\alpha\Vert v_i - z \Vert^2 + s_{bj}^\alpha\Vert v_j - z \Vert^2,
\end{equation}
which minimized with respect to $z$ ($\partial F/ \partial z = 0$), and assuming that $s_{kb}=s_{ib}+s_{jb}$, gives the following location of the branching point:
\begin{equation}
z_b = \frac{s_{bi}^\alpha v_i + s_{bj}^\alpha v_j + (s_{bi} + s_{bj})^\alpha v_k}{s_{bi}^\alpha + s_{bj}^\alpha + (s_{bi} + s_{bj})^\alpha}.
\end{equation}
Therefore, in our fluid network problem we should consider the mapping of the distribution of sectional areas of sources to the distribution of sectional areas of targets $\mu(s) \rightarrow \nu(s)$. 
We should mention also that the sectional areas are related to the volumetric flow rate: $q=dV/dt=sd\ell /dt= su$, where $u$ is the velocity of fluid, which are all measurable quantities.

Le us now consider a different approach to branching, by defining the following functional:
\begin{equation}
F(z) = (1-\alpha)(s_{bi}\Vert v_i - z \Vert^2 + s_{bj}\Vert v_j - z \Vert^2) + \alpha s_{kb}\Vert v_k - z \Vert^2,
\end{equation}
where $\alpha \in [0,1]$.  
Assuming that $s_{kb}=s_{bi}+s_{bj}$, the minimization gives the following solution for the branching point:
\begin{equation}
z_b = (1- \alpha )\frac{s_{bi} v_i + s_{bj} v_j}{s_{bi} + s_{bj}} + \alpha v_k.
\end{equation}
One can easily see that the resulted Y branching solution (17) is actually interpolating between a T like branching solution, obtained for $\alpha=0$, 
and a V like solution, corresponding to the OT problem, and obtained for $\alpha = 1$ (Figure 1).

\begin{figure}[!h]
\centering \includegraphics[width=12cm]{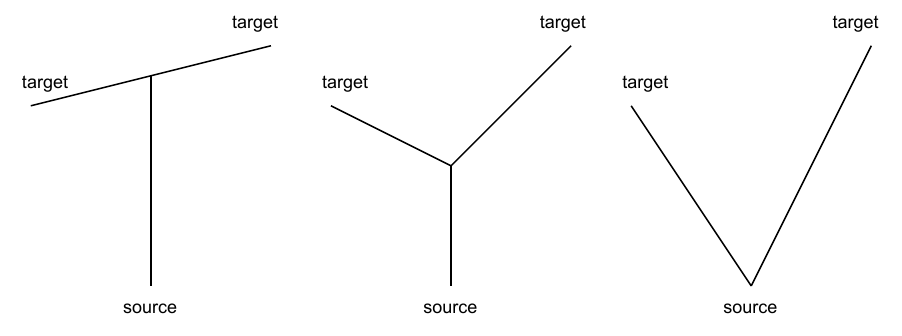}
\caption{The Y branching solution interpolating between T and V.}
\end{figure}

\section{Network optimization} 

We approach the network optimization in two main steps. In the first step we solve the OT problem (6)-(7) using the linear programming approach, or a variant of the Sinkhorn algorithm. 
This means that we obtain an optimal set of targets associated with each source, and we also satisfy the sectional area constraints. In the second step we use the local optimization approach described in the previous section 
to further minimize the cost of mapping each source to its associated targets. 
This way the problem is decomposed into $N_x$ one-to-many problems, where $N_x$ is the number of sources. 

In what follows we will formulate a heuristic algorithm to approximately solve the 
one-to-many problem in a deterministic number of steps. 
Let us assume that the source is at location $v_0$ and it has $N$ allocated targets: $V=[v_n]_{n=1}^N$. The goal is to find the optimal set of 
branching points with locations $\{v_m\}_{m=N+1}^{N+K}$ in order to minimize the flow from the source to the targets. We should determine the location of the branching points, however their number and sectional areas are also unknown. 
With each node we also associate a variable $a_n\in \{0,1\}$, where $a_n=1$ if the node $n$ is selectable, and $a_n=0$ otherwise. 
Also, with each node $n$ we associate a variable $c_n=m$ if node $m$ is connected to node $n$. After the initial OT step all the targets are connected to the source, so initially $c_n=0$ for $n=0,...,N$. 
We consider the following algorithm to find the branching points:
\begin{enumerate}
\item Initialize the network: $V=[v_n]_{n=0}^N$, $A=[a_n=1]_{n=0}^N$, $C=[c_n=0]_{n=0}^N$, $S=[s_{0n}]_{n=0}^N$, $\alpha$.
\item If $\sum_{n=1}^N a_n = 0$ return $V$ and $S$.
\item Find the node $i$ with $i=\text{arg}\max_i \Vert v_i - v_0\Vert$ and $a_i=1$. 
\item Find the closest node $j\neq i$ to $i$ for which the branching point $z_b$ computed with (17) (or 15) satisfies:
$j=\text{arg}\min_{j\neq i} \Vert v_j - v_i\Vert$ and $a_j=1$.
\item If $j>0$ then: $v_{N+1}=z$, $a_{N+1}=1$, $a_i=0$, $a_j=0$, $c_{i}=N+1$, $c_j=N+1$, $c_{N+1}=0$, $s_{0,N+1}=s_{0i}+s_{0j}$.  
\item If $j=0$ set $a_i=0$, $c_i=0$.
\item Go to step 2.
\end{enumerate}

At each step, the second stage of the algorithm proceeds by finding the optimal branching point for the farthest selectable target and its closer neighbor. Once these are found a branching point is computed and added to the network, 
and the targets are set to unselectable. 
From this point of view this is a combination of greedy and tabu search \cite{key-13}, since previously selected nodes cannot be re-selected again. 
The algorithm will end in a finite number of steps, when no other branching node can be computed. Each node $v_n$ with $n>N$ is a branching node in the returned list $V=[v_n]_{n=0}^{N+K}$. Branching can also 
occur between a branching node and a target, or two branching nodes. In the end all the targets will be connected to the source directly, or via some branching nodes. 

The algorithm's complexity is $O(N^2)$, since it requires pairs of points, and it has the advantage that it is deterministic and finite, and with each branching node added to the network the overall cost is guaranteed to decrease,   
showing very good results in practical simulations. 

\section{Numerical examples}

\subsection{Single source and multiple targets}

Let us first consider the case of a source situated at $v_0=[0,0]$ and $N$ targets with the positions randomly generated in $[-1,1]^3$. 
The sectional areas $s_{0n}>0$, $n=1,...,N$ are also randomly set in $[0,1]$ such that $\sum_n s_{0n}=1$. In Figure 2 (left) we give an example 
for $N=100$ targets and $\alpha=0.5$, using the local branching optimization using the equation (17). 
The source is represented by the larger red dot in the center, and the targets are the blue dots. 
The result for the transportation cost minimization is also shown in Figure 2 (right) as a function of the iteration step. One can see that the cost is monotonously 
decreasing, and in less than 100 steps it converges to a quite reasonable BOT solution.  

\begin{figure}[!h]
\centering \includegraphics[width=7.5cm]{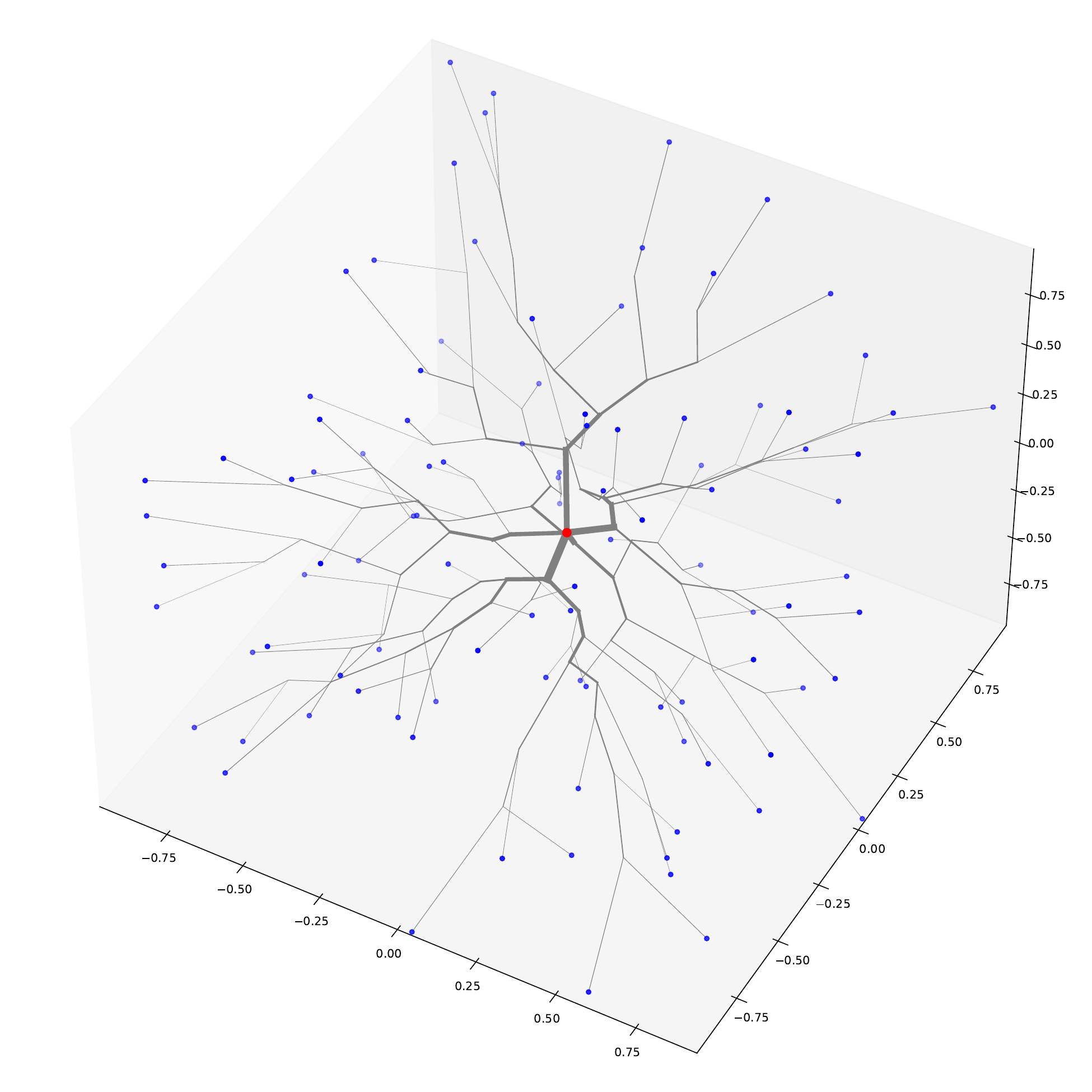}
\centering \includegraphics[width=7.5cm]{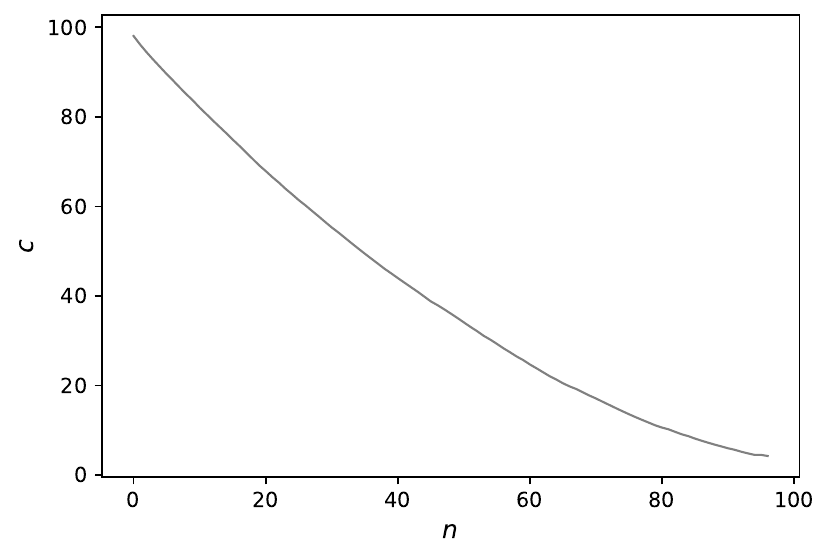}
\caption{Transport network and cost for a single source and N random targets.}
\end{figure}

\subsection{A simplified "cardiovascular" network}

In the second example we consider a simplified "cardiovascular" network, where the source corresponds to the heart, and we have two transport networks, one for arteries and the other for veins. 
Here we are simulating a closed circuit, where the blood flows from the heart to the targets (organs, muscles etc.) on the arteries, and then it returns back from the targets to the heart 
on the veins network. The two networks should be different, in the sense that there is some symmetry at the targets level, but then it is gradually lost since the veins and arteries must be separated in space. 

In order to model the artery-vein separation we add a "shifting" term to the equation (17):
\begin{equation}
z = (1 -\alpha) \frac{s_{bi} v_i + s_{bj} v_j}{s_{bi} + s_{bj}} + \alpha v_k + \frac{\varepsilon}{s_{bi} + s_{bj} + \delta}.
\end{equation}
Here, $\varepsilon$ is a random vector with a small norm $0<\Vert \varepsilon \Vert \ll 1$, and it has the same dimensionality as the nodes. Also, $0<\delta \ll 1$ is a small constant. 

\begin{figure}[!ht]
\centering \includegraphics[width=7.5cm]{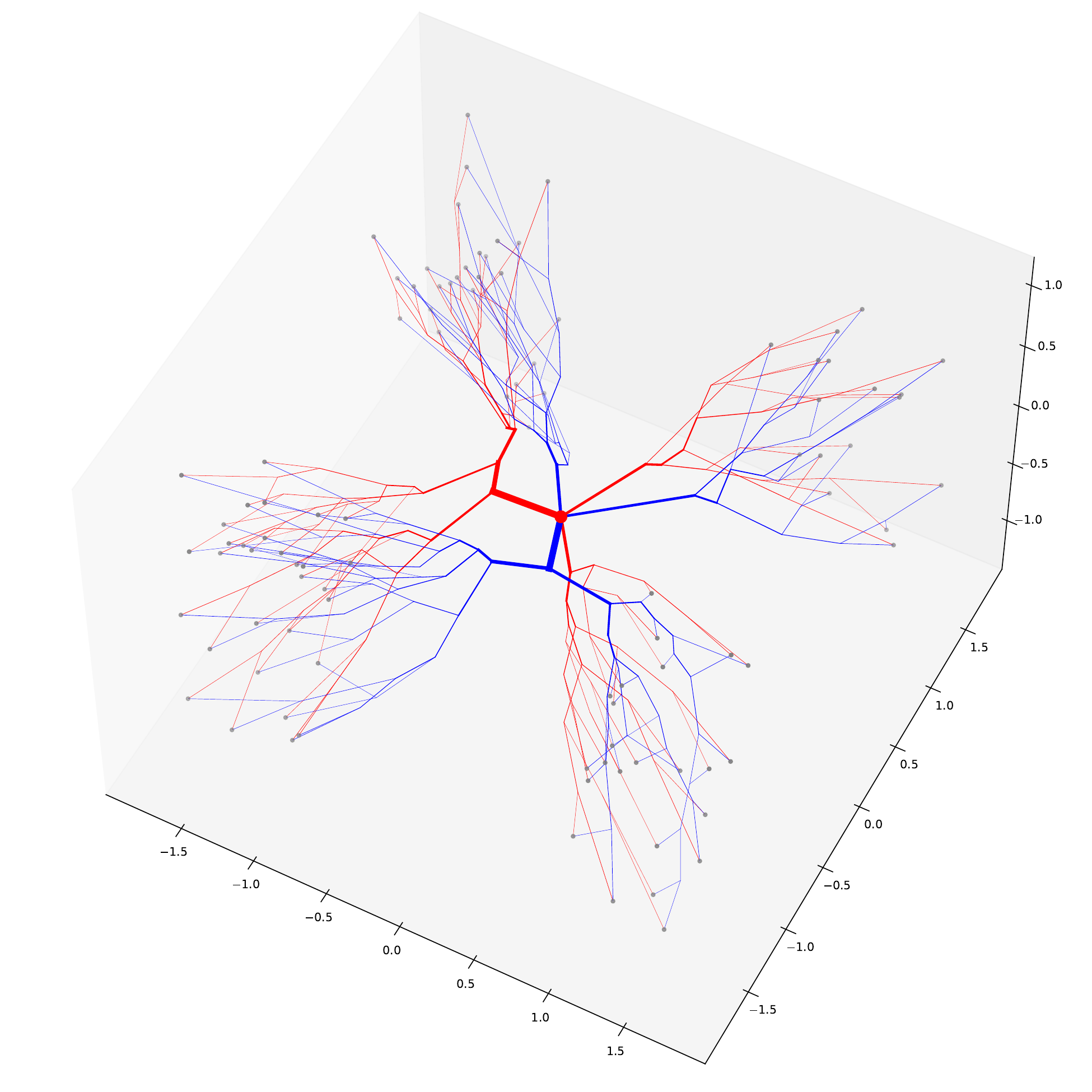}
\centering \includegraphics[width=7.5cm]{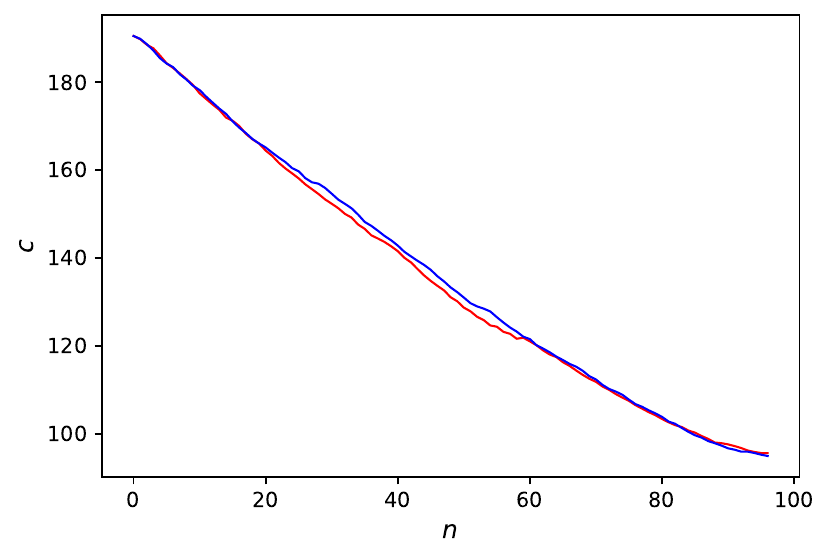}
\caption{Transport network and cost for the cardiovascular network (arteries=red, veins=blue).}
\end{figure}

The shift vector $\varepsilon$ in equation (18) is generated only once, at the beginning of the computation, and then it is kept frozen. The reason for the inverse proportionality of the shifting term with 
the quantity $s_{bi} + s_{bj}$ is that blood vessels with a larger sectional area are more rigid, while blood vessels with small sectional areas should be more susceptible to directional shifting and bending. 
The simulation results shown in Figure 3 have been obtained with the following parameters: $\alpha=0.5$, and $\Vert \varepsilon \Vert = \delta = 0.01$. 

\begin{figure}[!ht]
\centering \includegraphics[width=7.5cm]{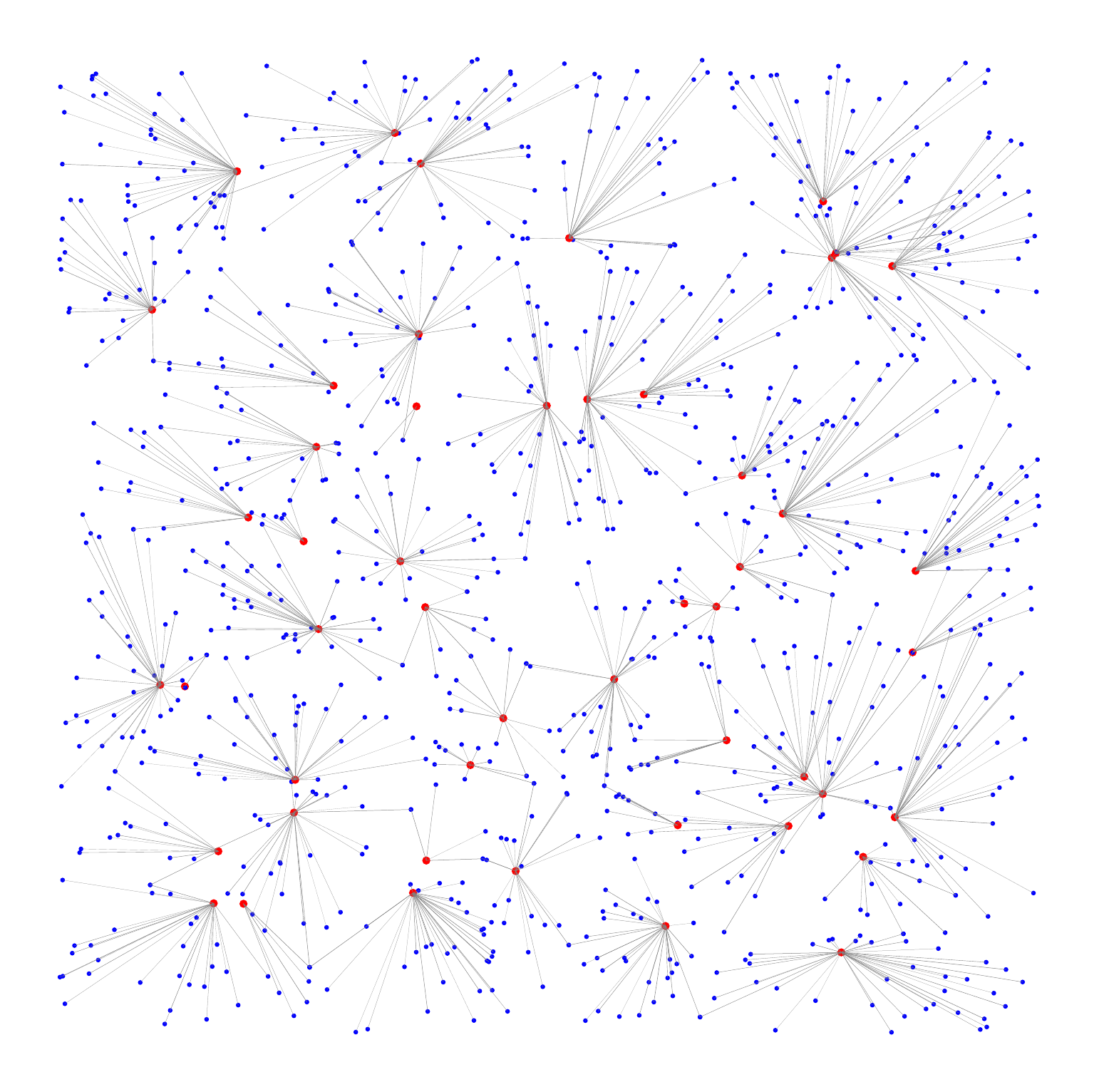}
\centering \includegraphics[width=7.5cm]{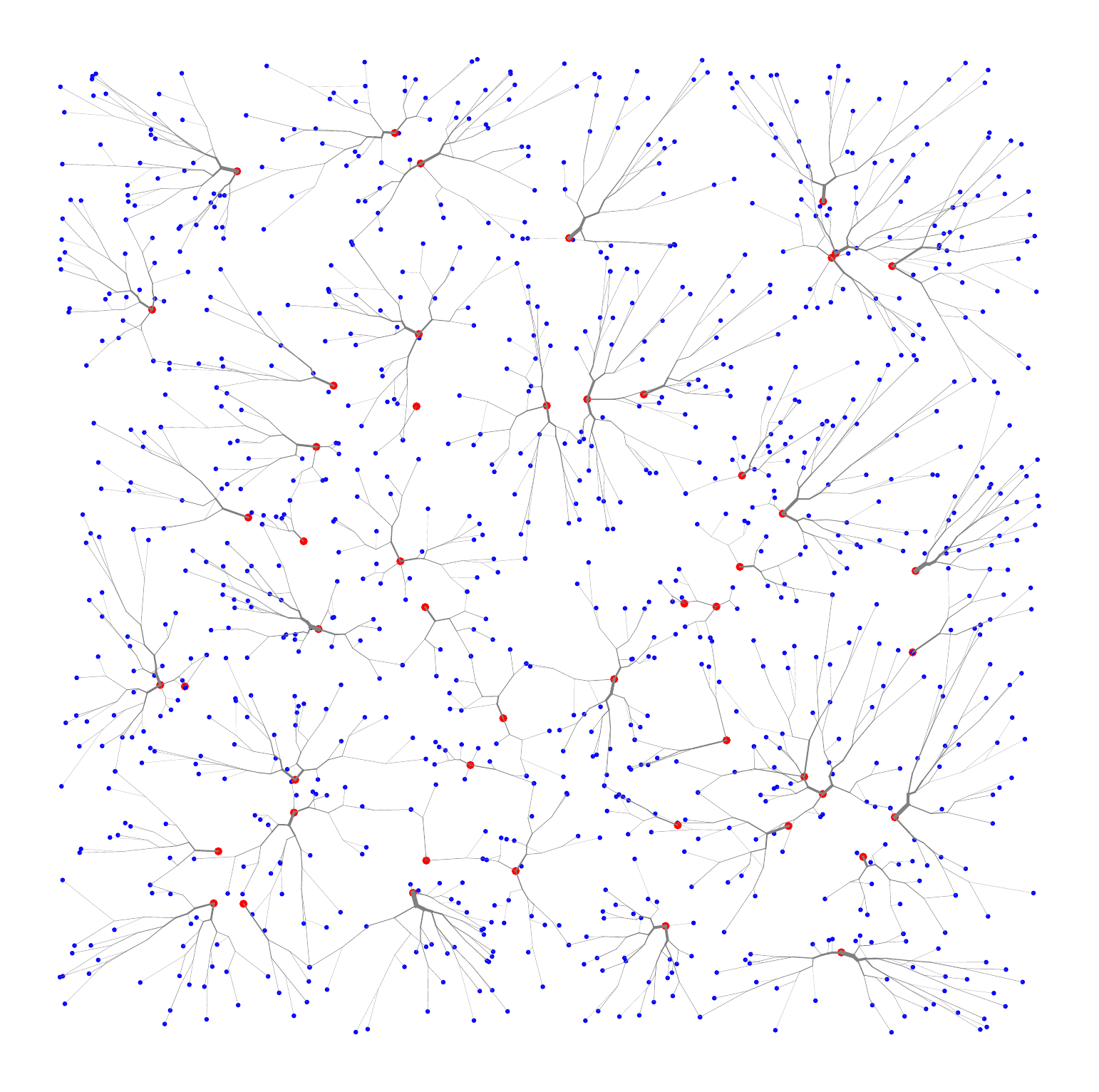}
\caption{Multiple sources and multiple targets: (left) the linear programming OT solution; (right) the BOT solution (sources=red, targets=blue).}
\end{figure}

\subsection{Multiple sources and multiple targets}

Let us now consider a more general case with multiple sources and targets, where $N_x$ sources and $N_y$ targets are randomly generated with positions in $[-1,1]^2$. 
The sectional areas are also randomly generated, and they satisfy the imposed constraints (3) and (4). 

In Figure 4 we give an example 
for $N_x=50$ sources and $Ny=1000$ targets, and $\alpha=0.25$. On the left of Figure 5 we give the linear programming solution of the OT problem, and on the right of Figure 5 we give the 
BOT solution. The sources are the red dots, and the targets are the blue dots. In this case, the cost was reduced from $c=5.11$, for the OT solution, to $c=0.63$ for the BOT solution. 

\subsection{Santa's transportation network}

Let us now consider Santa's transportation problem. Every Christmas, Santa Claus must transport and deliver the gifts to the people in each city around the World, starting from a North Pole undisclosed location:) 

In this problem we are given 141,182 cities around the world, with known location and population. 
In order to simplify the problem, we consider a national distribution center for each country which is at the center of (population) mass of the country. This can be easily computed since 
we have the position and the population of the cities:
\begin{equation}
x^*_0 = \frac{\sum_i p_i x_i}{\sum_i p_i},
\end{equation}
here $x_i$ is the location of the city and $p_i$ is the population. 
In addition we also define several regional distribution centers $\{x^*_1,x^*_2,...,x^*_K\}$. The number of centroids is empirically set to $K=\sqrt{N}+1$, where $N$ is the number of cities in the country, 
and their position is given by the centroids solution of the weighted K-means clustering:
\begin{equation}
\text{arg}\min_{\{x^*_1,x^*_2,...,x^*_K\}} \sum_k \sum_i p_i\Vert x_i - x^*_k\Vert^2,
\end{equation}
where the centroids $x^*_k$ are the center of mass of the clusters $C_k$:
\begin{equation}
x^*_k = \frac{\sum_i p_i x_i}{\sum_i p_i}, \quad i \in C_k.
\end{equation}
After we have computed the location of the national and regional distribution centers, we use the BOT algorithm to compute the branching points and to find a sub-optimal distribution network. 
The results are quite good as one can see in the web application demonstration at: \url{https://mandrecut.github.io/santa_net}. 
For example, in Figure 5 we show the global distribution network, connecting the North Pole with each country's centroid, and in Figure 6 we give an example for the United States distribution network, 
computed with the BOT algorithm. 
We also make available the data used in this example, since it could be used as a test case for other algorithms, 
and future developments (\url{https://github.com/mandrecut/santa_net}). 

\begin{figure}[!ht]
\centering \includegraphics[width=15cm]{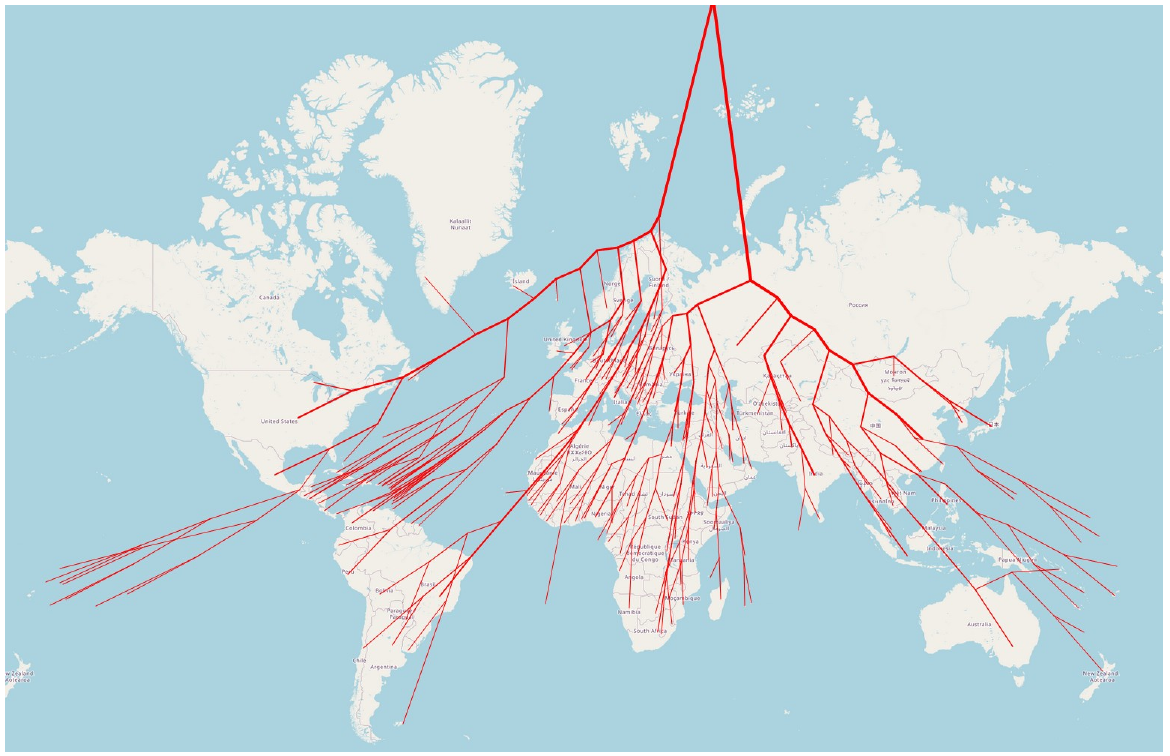}
\caption{Santa's global distribution network.}
\end{figure}

\begin{figure}[!ht]
\centering \includegraphics[width=15cm]{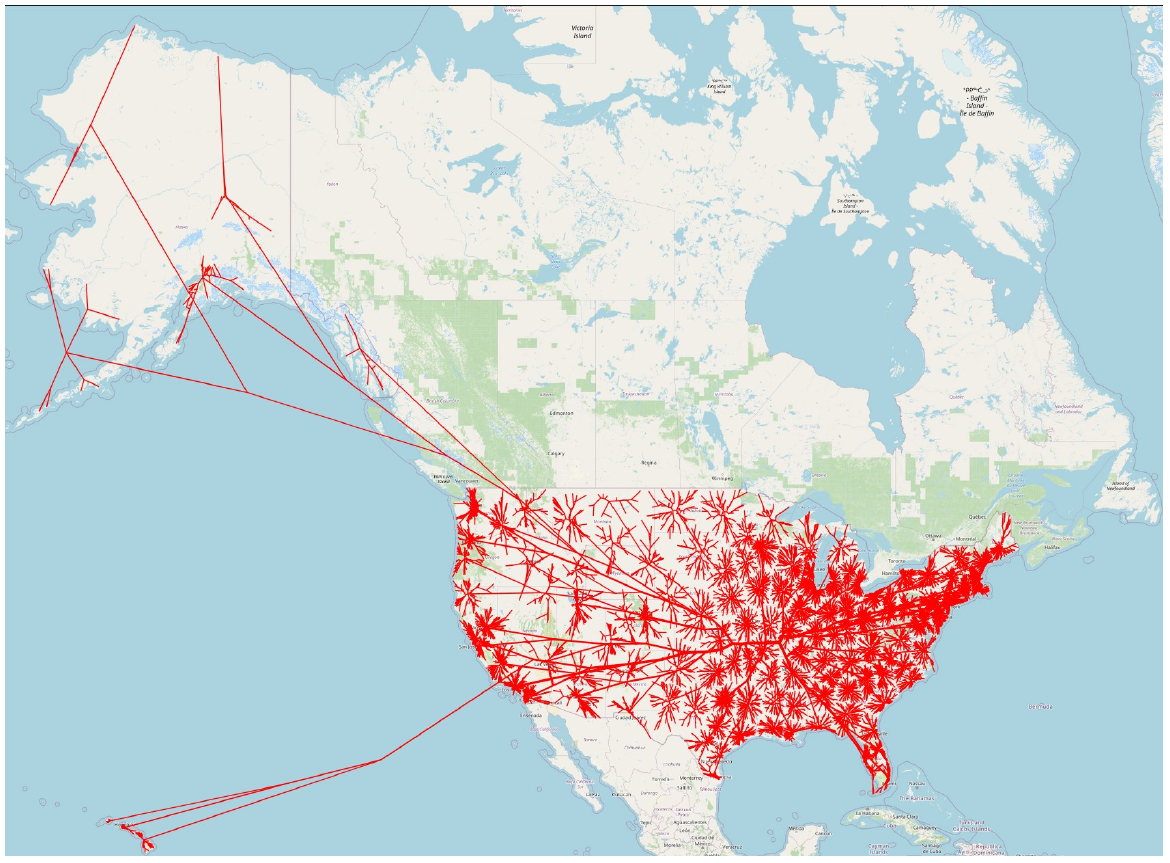}
\caption{Santa's distribution network for United States.}
\end{figure}

\section*{Conclusion}

We have described a heuristic branching method for optimal transport, inspired by fluid networks. The method is 
based on a deterministic algorithm implementing a combination of (regularized) linear programming in the first stage, followed by a greedy tabu search in the second stage, exploiting 
a local branching optimization approach. Due to its high speed 
and reduced complexity it can be successfully used to provide fast and good approximate solutions to larger scale problems. 
We have also discussed several numerical applications, including synthetic examples, a simplified cardiovascular network, and the "Santa Claus" distribution network which includes 141,182 cities around the world, 
with known location and population. The data and some of the code is available at: \url{https://github.com/mandrecut/santa_net}.

\end{document}